\documentclass[12pt,twoside]{article}

\usepackage{amssymb}
\usepackage{amsmath}
\usepackage{bbm}
\usepackage{mathrsfs}
\usepackage{float}
\usepackage{xypic}

\makeatletter
\def\hsmash{\relax 
  \ifmmode\def\next{\mathpalette\mathhsm@sh}\else\let\next\makehsm@sh
  \fi\next}
\def\makehsm@sh#1{\setbox\z@\hbox{#1}\finhsm@sh}
\def\mathhsm@sh#1#2{\setbox\z@\hbox{$\m@th#1{#2}$}\finhsm@sh}
\def\finhsm@sh{\wd\z@\z@ \box\z@}
\makeatother

\makeatletter
\gdef\th@mychange{\normalfont\slshape
   \def\@begintheorem##1##2{\item
        [\hskip\labelsep \theorem@headerfont ##2. ##1  \,--\!--\!--\!--  ]}%
 \def\@opargbegintheorem##1##2##3{%
   \item[\hskip\labelsep \theorem@headerfont ##2. ##1\ {\upshape(}##3{\upshape)}. \,-----  ]}}
\makeatother

\RequirePackage{theorem}
\theoremstyle{mychange}

{\theorembodyfont{\rmfamily}\newtheorem{ttt}{}[section]}
{\theorembodyfont{\rmfamily}\newtheorem{nota}[ttt]{Notation.}}
{\theorembodyfont{\rmfamily}\newtheorem{defi}[ttt]{Definition.}}
{\theorembodyfont{\rmfamily}\newtheorem{defis}[ttt]{Definitions.}}
{\theorembodyfont{\rmfamily}\newtheorem{remark}[ttt]{Remark.}}
{\theorembodyfont{\rmfamily}\newtheorem{rems}[ttt]{Remarks.}}
{\theorembodyfont{\rmfamily}\newtheorem{ex}[ttt]{Example.}}
{\theorembodyfont{\rmfamily}\newtheorem{fade}[ttt]{Fact-Definition.}}

{\theorembodyfont{\itshape}\newtheorem{fac}[ttt]{Fact.}}
{\theorembodyfont{\itshape}\newtheorem{lem}[ttt]{Lemma.}}
{\theorembodyfont{\itshape}\newtheorem{prop}[ttt]{Proposition.}}
{\theorembodyfont{\itshape}\newtheorem{coro}[ttt]{Corollary.}}

{\theorembodyfont{\rmfamily}\newtheorem{defio}[ttt]{Definition}}
{\theorembodyfont{\rmfamily}\newtheorem{algoo}[ttt]{Algorithm}}

{\theorembodyfont{\itshape}\newtheorem{propo}[ttt]{Proposition}}
{\theorembodyfont{\itshape}\newtheorem{theoo}[ttt]{Theorem}}

\newcounter{abc}
\newenvironment{abc}{\begin{list}{\rm \alph{abc}) }{\usecounter{abc} \leftmargin=0.0pt \labelsep=0.0pt \listparindent=0.0pt \labelwidth=0.0pt \parsep=\smallskipamount \itemsep=0.0pt \topsep=0.0pt \partopsep=\smallskipamount}}{\end{list}}
\newcounter{iii}
\newenvironment{iii}{\begin{list}{\rm \roman{iii}) }{\usecounter{iii} \leftmargin=0.0pt \labelsep=0.0pt \listparindent=0.0pt \labelwidth=0.0pt \parsep=\smallskipamount \itemsep=0.0pt \topsep=0.0pt \partopsep=\smallskipamount}}{\end{list}}

\newcommand{\bP}{\mathop{\text{\bf P}}\nolimits}
\newcommand{\Spec}{\mathop{\text{\rm Spec}}\nolimits}
\newcommand{\tr}{\mathop{\text{\rm tr}}\nolimits}
\newcommand{\N}{\mathop{\text{\rm N}}\nolimits}
\newcommand{\im}{\mathop{\text{\rm im}}\nolimits}
\newcommand{\Gal}{\mathop{\text{\rm Gal}}\nolimits}
\newcommand{\NT}{\mathop{\text{\rm NT}}\nolimits}
\newcommand{\Res}{\mathop{\text{\rm Res}}\nolimits}
\newcommand{\disc}{\mathop{\text{\rm disc}}\nolimits}
\newcommand{\pr}{\mathop{\text{\rm pr}}\nolimits}
\newcommand{\Pic}{\mathop{\text{\rm Pic}}\nolimits}
\newcommand{\Br}{\mathop{\text{\rm Br}}\nolimits}
\renewcommand{\L}{\mathop{\text{\rm L}}\nolimits}
\newcommand{\cha}{\mathop{\text{\rm char}}\nolimits}

\newcommand{\calO}{\mathscr{O}}

\newcommand{\bbC}{{\mathbbm C}}
\newcommand{\bbL}{{\mathbbm L}}
\newcommand{\bbQ}{{\mathbbm Q}}
\newcommand{\bbZ}{{\mathbbm Z}}

\newcommand{\id}{{\rm id}}

\newcommand{\br}{ }
\newcommand{\brr}{, }

\def\rightend#1#2{{%
 \leavevmode\nobreak\hskip .5em plus 1fil
 \penalty600 \hskip 0pt plus -1filll
 \vadjust{}\nobreak\hskip 0pt plus 1filll%
 #1\parfillskip=#2\relax \par}}

\def\eop{\ifmmode\rule[-22pt]{0pt}{1pt}\ifinner\tag*{$\square$}\else\eqno{\square}\fi\else\rightend{$\square$}{0pt}\fi}

\newcommand{\ratarrow}{$%
$\definemorphism{rat}\dashed\tip\notip%
\spreaddiagramcolumns{-12pt}%
\! - \!\!\diagram%
\rrat & 
\enddiagram\!\!$%
$}

\renewcommand{\thefootnote}{\arabic{footnote}}
\author{Andreas-Stephan Elsenhans${}^*$ and J\"org Jahnel${}^*$}
\date{}

\title{Cubic surfaces with a Galois invariant pair of~Steiner~trihedra}

\begin{document}
\renewcommand{\thefootnote}{\fnsymbol{footnote}}

\maketitle

\begin{abstract}
We present a method to construct non-singular cubic surfaces
over~$\bbQ$
with a Galois invariant pair of Steiner~trihedra. We~start with cubic surfaces in a form generalizing that of A.~Cayley and G.~Salmon. For~these, we develop an explicit version of Galois~descent.
\end{abstract}


\footnotetext[1]{
The computer part of this work was executed on the Sun Fire V20z Servers of the Gau\ss\ Laboratory for Scientific Computing at the G\"ottingen Mathematical Institute. Both authors are grateful to Prof.~Y.~Tschinkel for the permission to use these machines as well as to the system administrators for their~support.}

\section{Introduction}

\begin{ttt}
The~configuration of the 27~lines upon a smooth cubic surface is highly~symmetric. The~group of all permutations respecting the intersection pairing is isomorphic to the Weyl
group~$W(E_6)$
of
order~$51\,840$.

When~$S$
is a cubic surface
over~$\bbQ$,
the absolute Galois
group~$\Gal(\overline\bbQ/\bbQ)$
operates on the 27~lines. This~yields a
subgroup~$G \subseteq W(E_6)$.
\end{ttt}

\begin{ttt}
There~are exactly 350 conjugacy classes of subgroups
of~$W(E_6)$.
Only~for about one half of them, explicit examples of cubic surfaces
over~$\bbQ$
are~known.\smallskip

General~cubic surfaces~\cite{EJ1} lead to the
full~$W(E_6)$.
In~\cite{EJ2}, we constructed examples for the index two subgroup which is the simple group of
order~$25\,920$.
Other~examples may be obtained by fixing a
$\bbQ$-rational
line or tritangent~plane. Generically,~this yields the maximal subgroups
in~$W(E_6)$
of indices 27 and 45,~respectively. It~is not yet clear which smaller groups arise by further~specialization.

On~the other hand, there are a number of rather small subgroups
in~$W(E_6)$
for which examples may be constructed easily. Blowing~up six points
in~$\bP^2_\bbQ$
forming a Galois invariant set leads to a cubic surface with a Galois invariant~sixer. It~is clear that examples for all the 56 corresponding conjugacy classes of subgroups may be constructed in this~way. There~are a few more trivial cases, e.g.~diagonal surfaces, but all in all not more than 70 of the 350 conjugacy classes of subgroups may be realized by such elementary~methods.\smallskip

In~\cite{EJ3}, we presented a method to construct cubic surfaces
over~$\bbQ$
with a Galois invariant double-six. A~simple calculation in~{\tt GAP} shows that there are 102~conjugacy classes of subgroups
of~$W(E_6)$
fixing a double-six but no~sixer. For~each of them, explicit~examples of cubic surfaces are given in the list~\cite{EJ3a}. The~most interesting ones were reproduced in~\cite{EJ3}.
\end{ttt}

\begin{ttt}
In~this article, we present a method to construct cubic surfaces
over~$\bbQ$
with a Galois invariant pair of Steiner~trihedra.
Our~method is based on cubic~surfaces in a form slightly generalizing that of A.~Cayley and G.~Salmon. For~these, we develop an explicit version of Galois~descent.

There~are 63~conjugacy classes of subgroups
of~$W(E_6)$
which fix a pair of Steiner~trihedra but no~double-six. We~constructed explicit examples of cubic surfaces for each of~them. Some~of the most interesting ones will be presented in the final~section.
\end{ttt}

\section{Steiner trihedra}

This~section will mainly recall definitions and facts which are necessary for the~sequel. Most~of them were known to the geometers of the 19th~century~\cite{St,Do}.

\begin{ttt}
Let~$S$
be a smooth cubic surface over an algebraically closed~field. It~is well-known that
$S$
contains a total of 27~lines. There~are exactly 45~planes cutting three lines out
of~$S$.
These~are called the {\em tritangent~planes}.

Two~planes
in~$\bP^3$
which are different from each other always meet in a single~line. Given~two tritangent planes, there are two~possibilities. Their~intersection is either one of the 27~lines contained
in~$S$
or a line not contained
in~$S$.
For~a tritangent
plane~$E$,
there are twelve tritangent planes meeting
$E$
within the~surface, four for each of the lines
in~$E \cap S$.
32~tritangent planes remain which meet
$E$~otherwise.
\end{ttt}

\begin{remark}
The~set of pairs of distinct tritangent planes is acted upon by the automorphism
group~$W(E_6)$.
Under~this operation, that set is decomposed into exactly two orbits according to the way the tritangent planes meet each~other.
\end{remark}

\begin{defis}
\begin{abc}
\item
A~{\em trihedron\/} consists of three distinct tritangent planes such that the intersection of any two is not contained
in~$S$.
\item
For a trihedron
$\{ E_1, E_2, E_3 \}$,
a plane
$E$
is called a {\em conjugate plane\/} if each of the lines
$E_1 \cap E$,
$E_2 \cap E$,
and~$E_3 \cap E$
is contained in the
surface~$S$.
\end{abc}
\end{defis}

\begin{fade}
A~trihedron may have either no, exactly one, or exactly three conjugate~planes.
Correspondingly, a trihedron is said to be of the {\em first kind,} {\em second kind,} or {\em third~kind}. Trihedra~of the third kind are also called {\em Steiner~trihedra}.
\end{fade}

\begin{remark}
Let two tritangent planes
$E_1, E_2$
be given such that their intersection line is not contained in the
surface~$S$.
Then,~there are three tritangent planes meeting both
$E_1$
and~$E_2$
in lines
within~$S$.
Nine~further tritangent planes meet
$E_1$
on~$S$.
Analogously,~nine tritangent planes only meet
$E_2$
within~$S$.

22~tritangent planes~remain. Twelve~of them complete
$\{E_1, E_2\}$
to a trihedron of the first~kind. For~nine tritangent
planes~$E$,
$\{E_1, E_2, E\}$
becomes a trihedron of the second~kind. Finally,~there is a unique tritangent plane such that
$\{E_1, E_2, E\}$
is a Steiner~trihedron.

Consequently,~on a smooth cubic surface, there are 2880 trihedra of the first kind, 2160 trihedra of the second~kind, and 240 Steiner~trihedra. The~group
$W(E_6)$
acts transitively on the set of all Steiner~trihedra. In~fact, the operations on trihedra of the first and second kinds are transitive,~too.
\end{remark}

\begin{fac}
\label{compl}
\begin{abc}
\item
Steiner~trihedra come in~pairs. Actually,~the three conjugate planes of a Steiner trihedron form another Steiner~trihedron.
\item
Two trihedra define the same sets of lines if and only if they form a pair of Steiner~trihedra.
\item
The~nine lines defined by a Steiner trihedron form the complement of the lines contained in a triple of azygetic double-sixes.
\end{abc}\smallskip

\noindent
{\bf Proof.}
{\em
Recall~that two double-sixes on a non-singular cubic surface may be either syzygetic or azygetic according to the number of lines they have in~common. Further,~a pair of azygetic double-sixes uniquely determines a third double-six, azygetic to both of~them~\cite{Do,EJ4}.

The~assertion itself may best be seen in the blown-up model. In~Schl\"afli's notation~\cite[p.~116]{Sch}, one of the Steiner trihedra is formed by the tritangent planes
$[c_{14},c_{25},c_{36}]$,
$[c_{15},c_{26},c_{34}]$,
and~$[c_{16},c_{24},c_{35}]$.
Indeed,~the three conjugate planes are given by
$[c_{14},c_{26},c_{35}]$,
$[c_{15},c_{24},c_{36}]$,
and~$[c_{16},c_{25},c_{34}]$.
Further,~the~``standard'' triple of azygetic double-sixes
$$
\left(
\begin{array}{cccccc}
a_1 & a_2 & a_3 & a_4 & a_5 & a_6 \\
b_1 & b_2 & b_3 & b_4 & b_5 & b_6 
\end{array}
\right) \! ,
\;\;
\left(
\begin{array}{cccccc}
a_1    & a_2    & a_3    & c_{56} & c_{46} & c_{45} \\
c_{23} & c_{13} & c_{12} & b_4    & b_5    & b_6 
\end{array}
\right) \! ,
\;\; {\rm and\;~}
\left(
\begin{array}{cccccc}
c_{23} & c_{13} & c_{12} & a_4    & a_5    & a_6    \\
b_1    & b_2    & b_3    & c_{56} & c_{46} & c_{45} 
\end{array}
\right)
$$
is exactly formed by the remaining~lines.
}
\eop
\end{fac}

\begin{nota}
Let~$l_1, \ldots l_9$
be the nine lines defined by a Steiner~trihedron. Then,~we will denote the corresponding pair of Steiner~trihedra by a rectangular symbol of the form
$$\left[
\begin{array}{ccc}
l_1 & l_2 & l_3 \\
l_4 & l_5 & l_6 \\
l_7 & l_8 & l_9
\end{array}
\right] \! .$$
The~planes of the trihedra contain the lines noticed in the rows and~columns.
\end{nota}

\begin{prop}
\label{diag}
Let~a triple of azygetic double-sixes be~given.

\begin{abc}
\item
Then,~of the corresponding six sixers, one may form fifteen~pairs. Nine~of them are~disjoint. The~other six intersections are mutually disjoint triplets\/
$D_0, \ldots, D_5$.
\item
These~may be arranged in a diagram of the form
$$
\begin{array}{ccccc}
    & D_2 & \; & D_1 &  \\[1.5mm]
D_3 &     & &     & D_0 \\[1.5mm]
    & D_4 & & D_5
\end{array}
$$
such that two~triplets combine to a sixer if and only if they are~adjacent. Such~a diagram is unique up to rotation and~reflection.
\item
Further,~the following properties may be read off the~diagram.
\begin{iii}
\item
Every~line
in\/~$D_i$
meets every line
in\/~$D_j$
if and only if\/
$D_i$
and\/~$D_j$
are~opposite.
\item
Two~sixers form a double-six if and only if they are~opposite.
\item
The~nine lines in\/
$D_0 \cup D_2 \cup D_4$
are defined by a Steiner~trihedron. Analogously,~for the nine lines
in\/~$D_1 \cup D_3 \cup D_5$.
\end{iii}
\end{abc}

\noindent
{\bf Proof.}
{\em
Again,~let us work in the blown-up model and consider the standard triple of azgetic double-sixes formed by the 18~lines
$a_1, \ldots, a_6$,
$b_1, \ldots, b_6$,
$c_{12}$,
$c_{13}$,
$c_{23}$,
$c_{45}$,
$c_{46}$,
and~$c_{56}$.\smallskip

\noindent
Then,~a) is immediately~verified. The~six triplets which appear as intersections of the sixers are
$A_l := \{a_1, a_2, a_3\}$,
$A_r := \{a_4, a_5, a_6\}$,
$B_l := \{b_1, b_2, b_3\}$,
$B_r := \{b_4, b_5, b_6\}$,
$C_l := \{c_{12}, c_{13}, c_{23}\}$,
and~$C_r := \{c_{45}, c_{46}, c_{56}\}$.\smallskip

\noindent
b)
Consider~the~diagram
$$
\begin{array}{ccccc}
    & A_l & \; & A_r &  \\[1.5mm]
C_r &     & &     & C_l \, . \\[1.5mm]
    & B_l & & B_r 
\end{array}
$$
The~property stated is directly~checked. Uniqueness~is~clear.\smallskip

\noindent
c)
Properties~i) and~ii) may be verified~immediately. Further,~we have the two~pairs
$$\left[
\begin{array}{ccc}
a_1 & b_2 & c_{12} \\
b_3 & c_{23} & a_2 \\
c_{13} & a_3 & b_1
\end{array}
\right]
\qquad{\rm and}\qquad
\left[
\begin{array}{ccc}
a_4 & b_5 & c_{45} \\
b_6 & c_{56} & a_5 \\
c_{46} & a_6 & b_4
\end{array}
\right]$$
of Steiner~trihedra.
\eop
}
\end{prop}

\begin{fac}
Given a pair of Steiner trihedra, there is a unique way to decompose the 18~remaining lines into two sets of nine such that both are defined by Steiner~trihedra.\smallskip

\noindent
{\bf Proof.}
{\em 
The~{\em existence\/} of a decomposition as desired follows from Fact~\ref{compl}.c) and Proposition~\ref{diag}.c.iii). To~see {\em uniqueness}, we need an overview over all 120 pairs of Steiner~trihedra. In~the blown-up model, these are of the~types
$$
\left[
\begin{array}{ccc}
a_i & b_j & c_{ij} \\
b_k & c_{jk} & a_j \\
c_{ik} & a_k & b_i
\end{array}
\right] \! ,
\qquad
\left[
\begin{array}{ccc}
c_{il} & c_{jm} & c_{kn} \\
c_{jn} & c_{kl} & c_{im} \\
c_{km} & c_{in} & c_{jl}
\end{array}
\right] \! ,
\quad
{\rm and}
\quad
\left[
\begin{array}{ccc}
a_i & b_j & c_{ij} \\
b_k & a_l & c_{kl} \\
c_{ik} & c_{jl} & c_{mn}
\end{array}
\right] \! .
$$
We~have 20 pairs of Steiner~trihedra of the first type, 10 of the second, and 90 of the last~type. Having~seen this, it~is easy to verify that there are exactly two pairs of Steiner~trihedra having no line in common~with
$$\left[
\begin{array}{ccc}
c_{14} & c_{25} & c_{36} \\
c_{26} & c_{34} & c_{15} \\
c_{35} & c_{16} & c_{24}
\end{array}
\right] \! .\eop$$
}
\end{fac}

\begin{defi}
Given a pair of Steiner trihedra, we will call the two other pairs {\em complementary\/} to the given one if, altogether, they define all the 27~lines.
\end{defi}

\begin{rems}
\begin{iii}
\item
The~investigation above shows, in fact, that, for each pair of Steiner trihedra, there are exactly two pairs having no line in common, 54~pairs having two lines in common, 36~pairs having three lines in common, and 27~pairs which have five lines in~common with the nine lines defined by the pair~given.
\item
The~subgroup
of~$W(E_6)$
stabilizing a pair of Steiner trihedra is isomorphic
to~$[(S_3 \times S_3) \rtimes \bbZ/2\bbZ] \times S_3$
of
order~$432$.

A~calculation in {\tt GAP} shows that, indeed, this group operates on pairs of Steiner trihedra such that the orbits have lengths
$1$,~$2$,
$27$,
$36$,
and~$54$.
\end{iii}
\end{rems}

\section{The generalized Cayley-Salmon form}

\begin{nota}
One~way to write down a cubic surface explicitly is the so-called {\em Cayley-Salmon~form\/}~\cite[\S9.3]{Do}. A~slight generalization is the~following.
For~$u_0, u_1 \neq 0$,
denote
by~$\smash{S_{u_0, u_1}^{(a_0, \ldots, a_5, b_0, \ldots, b_5)}}$
the cubic surface given
in~$\bP^5$
by the system of~equations
\begin{eqnarray*}
u_0 X_0 X_1 X_2 + u_1 X_3 X_4 X_5 & = & 0 \,
, \\
a_0 X_0 + a_1 X_1 + a_2 X_2 + a_3 X_3 + a_4 X_4 + a_5 X_5 & = & 0 \, , \\
b_0 X_0 + \hsmash{b_1}\phantom{a_1} X_1 + \hsmash{b_2}\phantom{a_2} X_2 + \hsmash{b_3}\phantom{a_3} X_3 + \hsmash{b_4}\phantom{a_4} X_4 + \hsmash{b_5}\phantom{a_5} X_5 & = & 0 \, .
\end{eqnarray*}
\end{nota}

\begin{remark}
The~geometric meaning of these equations is to intersect the cubic fourfold given by
$u_0 X_0 X_1 X_2 + u_1 X_3 X_4 X_5 = 0$
with two hyperplanes. All~these fourfolds are actually isomorphic to each~other.
For~$u_0 = u_1 = 1$,
the classical Cayley-Salmon~form is~obtained.
\end{remark}

\begin{defi}
Let~$S_{u_0, u_1}^{(a_0, \ldots, a_5, b_0, \ldots, b_5)}$
be a cubic surface in generalized Cayley-Salmon form. We~will call the general cubic polynomial
\begin{eqnarray*}
 & & \Phi_{u_0, u_1}^{(a_0, \ldots, a_5, b_0, \ldots, b_5)}(T) := \\
 & & \hspace{2cm} \frac1{u_0} (a_0 + b_0 T)(a_1 + b_1 T)(a_2 + b_2 T) - \frac1{u_1} (a_3 + b_3 T)(a_4 + b_4T)(a_5 + b_5T)
\end{eqnarray*}
the {\em auxiliary polynomial\/} associated
with~$S_{u_0, u_1}^{(a_0, \ldots, a_5, b_0, \ldots, b_5)}$.
We~will simply write
$\Phi$
instead of~$\smash{\Phi_{u_0, u_1}^{(a_0, \ldots, a_5, b_0, \ldots, b_5)}}$
when there is no danger of~confusion.
\end{defi}

\begin{propo}[{\rm The discriminantal locus}{}]
\label{34}
Over~a base
field\/~$K$
of
characteristic\/~$\neq 3$,
the cubic surface\/
$S_{u_0, u_1}^{(a_0, \ldots, a_5, b_0, \ldots, b_5)}$
is singular if and only~if

\begin{iii}
\item
$$
\det
\left(
\begin{array}{cc}
a_i & a_j \\
b_i & b_j
\end{array}
\right)
= 0
$$
for some\/
$i \in \{0,1,2\}$
and\/~$j \in \{3,4,5\}$, or
\item
the discriminant of the auxiliary polynomial~vanishes.
\end{iii}\smallskip

\noindent
{\bf Proof.}
{\em
There~are two ways the intersection of the cubic fourfold given by
\begin{equation}
\label{CS}
u_0 X_0 X_1 X_2 + u_1 X_3 X_4 X_5 = 0
\end{equation}
with the two hyperplanes may become~singular. On~one hand, it might happen that both hyperplanes meet a singular point of the~fourfold.

The~singular locus of~(\ref{CS}) is given by
$$X_0X_1 = X_0X_2 = X_1X_2 = X_3X_4 = X_3X_5 = X_4X_5 = 0.$$
This means nothing but
$X_{i_1} = X_{i_2} = X_{j_1} = X_{j_2} = 0$
for
$i_1 \neq i_2 \in \{0,1,2\}$
and\/~$j_1 \neq j_2 \in \{3,4,5\}$.
We~meet such a point if and only if the corresponding determinantal condition is~fulfilled. The~degenerate case that the two linear forms are lineraly dependent is covered by this case,~too.

On~the other hand, the hyperplanes might meet the fourfold tangentially in a certain
point~$(x_0 : \ldots : x_5)$.
This~means that the tangent hyperplane of the fourfold
at~$(x_0 : \ldots : x_5)$
is a linear combination of the two hyperplanes~given. The~tangent hyperplane is given~by
$$u_0 (x_1x_2 X_0 + x_0x_2 X_1 + x_0x_1 X_2) + u_1 (x_4x_5 X_3 + x_3x_5 X_4 + x_3x_4 X_5) = 0 \, .$$
There~are two~cases.\medskip\pagebreak[3]

\noindent
{\em First Case:}
One of the coordinates
$x_0, \ldots, x_5$~vanishes.\smallskip

\noindent
Then,~in both summands
of~$(\ref{CS})$,
at least one factor must~vanish. Without~restriction,
suppose
$x_0 = x_3 = 0$.
The~tangent hyperplane is then given by
$M X_0 + N X_3 = 0$
for certain constants
$M$
and~$N$.
This~may be a linear combination of
$a_0 X_0 + a_1 X_1 + a_2 X_2 + a_3 X_3 + a_4 X_4 + a_5 X_5$
and
$b_0 X_0 + b_1 X_1 + b_2 X_2 + b_3 X_3 + b_4 X_4 + b_5 X_5$
only if
$\det ({a_2 \,a_5 \atop b_2 \,b_5}) = 0$.\medskip

\noindent
{\em Second Case:}
$x_0, \ldots, x_5 \neq 0$.\smallskip

\noindent
Let~the tangent hyperplane be given
by~$M_0 X_0 + \ldots + M_5 X_5 = 0$.
Then,~necessarily,~$\frac1{u_0} M_0M_1M_2 = \frac1{u_1} M_3M_4M_5$.
The~point of tangency~is
\begin{eqnarray*}
 & & (x_0 : \ldots : x_5) = \\
 & & \textstyle \hspace{1.2cm} \big( \frac1{u_0}M_1M_2 : \frac1{u_0}M_0M_2 : \frac1{u_0}M_0M_1 : (-\frac1{u_1}M_4M_5) : (-\frac1{u_1}M_3M_5) : (-\frac1{u_1}M_3M_4) \big) .
\end{eqnarray*}
We~suppose that the tangent hyperplane is a linear combination of the two linear forms~given.
Then,~$M_i = b_i$
or
$M_i = a_i + b_i t$,
for
some~$t$,
$i = 0, \ldots, 5$.
The~first variant may be interpreted
as~``$t = \infty$''.

The~conditions that
$(x_0 : \ldots : x_5)$
must be contained in both given hyperplanes may be rephrased~as
\begin{eqnarray*}
\textstyle \frac1{u_0} [a_0 (a_1+b_1t)(a_2+b_2t) + a_1 (a_0+b_0t)(a_2+b_2t) + a_2 (a_0+b_0t)(a_1+b_1t)] \hspace{0.95cm} & & \\
\textstyle {}- \frac1{u_1} [a_3 (a_4+b_4t)(a_5+b_5t) - a_4 (a_3+b_3t)(a_5+b_5t) - a_5 (a_3+b_3t)(a_4+b_4t)] & = & 0 \,\phantom{.}
\end{eqnarray*}
and
\begin{eqnarray*}
\textstyle \frac1{u_0} [b_0 (a_1+b_1t)(a_2+b_2t) + b_1 (a_0+b_0t)(a_2+b_2t) + b_2 (a_0+b_0t)(a_1+b_1t)] \hspace{0.95cm} & & \\
\textstyle {}- \frac1{u_1} [b_3 (a_4+b_4t)(a_5+b_5t) - b_4 (a_3+b_3t)(a_5+b_5t) - b_5 (a_3+b_3t)(a_4+b_4t)] & = & 0 \, .
\end{eqnarray*}
In~terms of the auxiliary polynomial, these two quadratic polynomials are
$3\Phi - t \Phi'$
and~$\Phi'$.
As~they have a common zero, we see that
$\Res_{2,2} (3\Phi - t \Phi', \Phi')$
must~vanish.

Let~us calculate this~resultant. First,~the leading coefficient
of~$\Phi'$
is equal
to~$3 (\frac1{u_0}b_0b_1b_2 - \frac1{u_1}b_3b_4b_5)$.
Hence,~according to the definition of the~resultant,
$$\Res_{2,2} (3\Phi - t \Phi', \Phi') = \frac{\Res_{3,2} (3\Phi - t \Phi', \Phi')}{3 (\frac1{u_0}b_0b_1b_2 - \frac1{u_1}b_3b_4b_5)} \, .$$
On~the other hand,~as
$t \Phi'$
is a multiple
of~$\Phi'$,
\begin{eqnarray*}
 & & \Res_{3,2} (3\Phi - t \Phi', \Phi') = \Res_{3,2} (3\Phi, \Phi') = \\
 & & \textstyle \hspace{5.5cm} 9 \Res_{3,2} (\Phi, \Phi') = -9 (\frac1{u_0}b_0b_1b_2 - \frac1{u_1}b_3b_4b_5) \disc(\Phi) \, .
\end{eqnarray*}
Consequently,~$\Res_{2,2} (3\Phi - t \Phi', \Phi') = -3 \disc(\Phi)$.

For~$u_0$
and
$u_1$
fixed, this is an irreducible polynomial in twelve~variables. The~corresponding component really occurs in the discriminantal variety as, for example,
\begin{eqnarray*}
u_0 x_0 + x_1 + x_2 + u_1 x_3 - x_4 - x_5 & = & 0 \\
                    2 u_1 x_3 + x_4 + x_5 & = & 0
\end{eqnarray*}
yields tangency
at~$(\frac1{u_0}:1:1:(-\frac1{u_1}):1:1)$
although we do not meet any of the nine determinantal~components.
}
\eop
\end{propo}

\begin{remark}
The~actual discriminant is a polynomial of
degree~$32$
in
$u_0$
and~$u_1$
and
bidegree~$(24,24)$
in the
$a$'s
and~$b$'s.
It~factors into the squares of the nine determinants
\smash{$\det ({a_i \,a_j \atop b_i \,b_j})$}
and~$u_0^{18} u_1^{18} \disc \Phi$.
The~necessity of taking the squares is motivated by~\cite[Theorem~2.12]{EJ2} and Corollary~\ref{gerade},~below.
\end{remark}

\section{Obvious and non-obvious lines}

\begin{ttt}
\label{41}
Let~$S = S_{u_0, u_1}^{(a_0, \ldots, a_5, b_0, \ldots, b_5)}$
be a non-singular cubic surface in generalized Cayley-Salmon~form.
Then,~on~$S$,
there are nine lines of the~type
$$L_{ij}: X_i = X_j = 0, \qquad i = 0,1,2, \;j = 3,4,5,$$
which we call the {\em obvious lines}.
\end{ttt}

\begin{fac}
The~linear forms
$X_0, \ldots, X_5$
define six tritangent planes
$E_0, \ldots , E_5$.
They~form a~pair
$$\big( \{E_0,E_1,E_2\}, \{E_3,E_4,E_5\} \big)$$
of Steiner~trihedra.
\end{fac}

\begin{remark}
The~generalized Cayley-Salmon~form therefore distinguishes one of the 120~pairs of Steiner~trihedra.
\end{remark}

\begin{remark}
The~situation here is analogous to that of a non-singular cubic surface in the hexahedral form of L.~Cremona~\cite{Cr} and Th.~Reye~\cite{Re}. A~cubic surface in hexahedral form is given
in~$\bP^5$
by a system of equations of the~type
\begin{eqnarray*}
\phantom{a_0} X_0^3 + \phantom{a_1} X_1^3 + \phantom{a_2} X_2^3 + \phantom{a_3} X_3^3 + \phantom{a_4} X_4^3 + \phantom{a_5} X_5^3 & = & 0 \,
, \\
\phantom{a_0}\hsmash{X_0}\phantom{X_0^3} + \phantom{a_1}\hsmash{X_1}
\phantom{X_1^3} + \phantom{a_2}\hsmash{X_2}\phantom{X_2^3} + \phantom{a_3}\hsmash{X_3}\phantom{X_3^3} + \phantom{a_4}\hsmash{X_4}\phantom{X_4^3} + \phantom{a_5}\hsmash{X_5}\phantom{X_5^3} & = & 0 \, , \\
a_0 \hsmash{X_0}\phantom{X_0^3} + a_1 \hsmash{X_1}\phantom{X_1^3} + a_2 \hsmash{X_2}\phantom{X_2^3} + a_3 \hsmash{X_3}\phantom{X_3^3} + a_4 \hsmash{X_4}\phantom{X_4^3} + a_5 \hsmash{X_5}\phantom{X_5^3} & = & 0 \, .
\end{eqnarray*}
Here,~there are the 15 obvious lines given~by
$X_{i_0} + X_{i_1} = X_{i_2} + X_{i_3} = 0$
for
$(i_0 i_1)(i_2 i_3)(i_4 i_5)$
a partition of the
set~$\{0, \ldots,5\}$.
\end{remark}

\begin{prop}
\label{43}
Let\/~$K$
be a field
s.t.~$\cha K \neq 3$,
$\smash{S = S_{u_0, u_1}^{(a_0, \ldots, a_5, b_0, \ldots, b_5)}}$
a non-singular cubic surface in generalized Cayley-Salmon~form
over\/~$K$,
and\/
$\lambda$
a zero of the auxiliary
polynomial\/~$\Phi$
associated
with\/~$S$.

\begin{abc}
\item
Then,~$S$
has a hexahedral form in the coordinates
\begin{eqnarray*}
Z_0 := -Y_0 + Y_1 + Y_2, \quad
Z_1 := Y_0 - Y_1 + Y_2, \quad
Z_2 := Y_0 + Y_1 - Y_2, \\
Z_3 := -Y_3 + Y_4 + Y_5, \quad
Z_4 := Y_3 - Y_4 + Y_5, \quad
Z_5 := Y_3 + Y_4 - Y_5 \hspace{1.5mm}
\end{eqnarray*}
for
$$Y_i := (a_i + b_i \lambda) X_i, \quad i = 0,\ldots,5\,.$$
\item
In~particular, six non-obvious lines
on\/~$S$
may be described~by
$$L^\lambda_\rho : Z_0 + Z_{\rho(0)} = Z_1 + Z_{\rho(1)} = 0$$
for\/~$\rho \colon \{0,1,2\} \to \{3,4,5\}$
a~bijection.
\end{abc}\medskip

\noindent
{\bf Proof.}
{\em
a)
We~have
$$\textstyle \frac1{u_0} (a_0 + b_0 \lambda)(a_1 + b_1 \lambda)(a_2 + b_2 \lambda) = \frac1{u_1} (a_3 + b_3 \lambda)(a_4 + b_4 \lambda)(a_5 + b_5 \lambda) \, .$$
Hence,~$S$
is given~by
\begin{eqnarray}
\label{eins}
Y_0Y_1Y_2 + Y_3Y_4Y_5 & = & 0 \,, \\
\label{zwei}
Y_0 + Y_1 + Y_2 + Y_3 + Y_4 + Y_5 & = & 0 \,,
\end{eqnarray}
and another linear~relation. (\ref{eins})~and~(\ref{zwei}) together~imply
\begin{eqnarray*}
(-Y_0 + Y_1 + Y_2)^3 + (Y_0 - Y_1 + Y_2)^3 + (Y_0 + Y_1 - Y_2)^3 \hspace{5cm} \\
{}+ (-Y_3 + Y_4 + Y_5)^3 + (Y_3 - Y_4 + Y_5)^3 + (Y_3 + Y_4 - Y_5)^3 = 0 \,.
\end{eqnarray*}
We~note that
$Z_0, \ldots, Z_5$
are projective coordinates, i.e., linearly~independent. For~that, the only point that requires attention is to verify
$a_i + b_i \lambda \neq 0$
for
all~$i$.
But,~as
$\lambda$
is a zero of the auxiliary polynomial, the opposite would imply
$a_i + b_i \lambda = 0$
for
some~$i \in \{0,1,2\}$
and
$a_j + b_j \lambda = 0$
for
some~$j \in \{3,4,5\}$.
Then,~$\smash{\det \big( {a_i \,a_j \atop b_i \, b_j} \big) = 0}$
and
$S$~were~singular.\smallskip

\noindent
b)
is~clear.
}
\eop
\end{prop}

\begin{prop}
\label{44}
Let\/~$K$
be a field
and\/
$S = S_{u_0, u_1}^{(a_0, \ldots, a_5, b_0, \ldots, b_5)}$
a non-singular cubic surface in generalized Cayley-Salmon~form
over\/~$K$.
(Suppose\/
${\rm char}\, K = 0$.)\smallskip

\noindent
For\/~$\lambda$
a zero of the auxiliary polynomial,
put\/~$\smash{\L^\lambda := \big\{ L^\lambda_\rho \mid \rho \colon \{0,1,2\} \stackrel\cong\longrightarrow \{3,4,5\} \big\}}$.
Further,~for\/
$e := \left\{ ({012\atop345}), ({012\atop453}), ({012\atop534})\right\}$
and\/
$o := \left\{ ({012\atop354}), ({012\atop435}), ({012\atop543})\right\}$,
$$\L^\lambda_e := \{ L^\lambda_\rho \mid \rho \in e \} \quad{\it and}\quad \L^\lambda_o := \{ L^\lambda_\rho \mid \rho \in o \} \, .$$

\begin{iii}
\item
Let~$\lambda_1 \neq \lambda_2$
be zeros of the auxiliary~polynomial. Then,~the sets\/ 
$\smash{\L^{\lambda_1}}$
and\/
$\smash{\L^{\lambda_2}}$
are~disjoint. In~particular, the nine obvious lines\/
$L_{ij}$
together with the 18~non-obvious lines\/
$L^{\lambda_i}_\rho$
form the set of all the 27~lines
on\/~$S$.
\item
Let\/~$\lambda_1, \lambda_2, \lambda_3$
be the zeroes of the auxiliary~polynomial. Then,~the~sets
$\L^{\lambda_i}_e$
and~$\L^{\lambda_i}_o$
are~triplets. In~the sense of~Proposition~\ref{diag}, they form the~diagram
$$
\begin{array}{ccccc}
    & \L^{\lambda_1}_e & \; & \L^{\lambda_2}_o &  \\[1.5mm]
\L^{\lambda_3}_o &     & &     & \L^{\lambda_3}_e \, . \\[1.5mm]
    & \L^{\lambda_2}_e & & \L^{\lambda_1}_o
\end{array}
$$
In~particular, the sets of~lines
$\L^{\lambda_1}_e \cup \L^{\lambda_2}_e \cup \L^{\lambda_3}_e$
and\/~$\L^{\lambda_1}_o \cup \L^{\lambda_2}_o \cup \L^{\lambda_3}_o$
are defined by the pairs of Steiner trihedra, complementary to the one~distinguished.
\end{iii}
\medskip

\noindent
{\bf Proof.}
{\em
i)
Assume first that the auxiliary polynomial defines
an~$S_3$-extension
of~$K$.
Then,~the lines
in~$L^{\lambda_1}$
are defined
over~$K(\lambda_1)$
and not
over~$K$.
Analogously,~the lines
in~$L^{\lambda_2}$
are defined
over~$K(\lambda_2)$,
not
over~$K$.
As~$K(\lambda_1) \cap K(\lambda_2) = K$,
this implies the~assertion.

Next,~suppose
that~$K = \bbC$.
The~family of all non-singular cubic surfaces in Cayley-Salmon form is defined over an open subscheme of
$$\Spec \bbC[A_0, \ldots, A_5, B_0, \ldots, B_5] = {\bf A}^{\!12}.$$
The~generic fiber is a surface defined
over~$K = \bbC(A_0, \ldots, A_5, B_0, \ldots, B_5)$.
It~is easy to check that the auxiliary polynomial is irreducible
over~$K$
and its discriminant is a non-square. Hence,~the assertion is true for the generic~fiber. Under~specialization, intersection numbers are~unchanged. In~particular, different lines can not specialize to the~same. The~assertion~follows.

For a general base
field~$K$,
we have that
$S^{(a_0, \ldots, a_5, b_0, \ldots, b_5)}$
is the base change of the corresponding surface over the finitely generated
field~$K' := \bbQ(a_0, \ldots, a_5, b_0, \ldots, b_5)$.
As~this field injects
into~$\bbC$,
the proof is~complete.\smallskip

\noindent
ii)
It~is readily checked that the sets described are indeed~triplets. I.e.,~that they consist of skew~lines. Further,~every line in
$\L^{\lambda_i}_e$
meets every line
in~$\L^{\lambda_i}_o$.
Thus,~these have to be placed in opposite~positions.

It~remains to be shown that
$\L^{\lambda_i}_e$
and~$\L^{\lambda_j}_e$
may not be~adjacent. For~this, again, we may assume that the auxiliary polynomial defines
an~$S_3$-extension
of~$K$.
Suppose~$\L^{\lambda_1}_e$
and~$\L^{\lambda_2}_e$
were~adjacent. Then,~together, they would form a~sixer. The~Galois operation ensures the same for
$\L^{\lambda_2}_e$
and~$\L^{\lambda_3}_e$
and, as well, for
$\L^{\lambda_3}_e$
and~$\L^{\lambda_1}_e$.
It~is, however, impossible that three entries in the diagram are pairwise~adjacent.
}
\eop
\end{prop}

\section{The norm-trace construction}

\begin{ttt}
Let~$R$
be a commutative ring with unit and
$A$
a commutative
$R$-algebra
which is \'etale and of finite~rank.
Then,~$A$
is, in particular, a locally free
$R$-module.
For an element
$a \in A$,
we have its norm and~trace. In~the free case, these are defined by
$\N (a) := \det_R (\cdot a \colon A \to A) \in R$
and~$\tr (a) := \tr_R (\cdot a \colon A \to A) \in R$.
The~general case is obtained by~gluing.

This~definition immediately generalizes to polynomials with coefficients
in~$A$.
In~fact,
$A[T_1, \ldots, T_n]$
is \'etale
over~$R[T_1, \ldots, T_n]$
of the same rank as
$A$
is
over~$R$.
\end{ttt}

\begin{defio}[{\rm The norm-trace construction}{}]
Let~$D$
be a commutative semisimple algebra of dimension two
over~$\bbQ$
and
$A$
a commutative
$D$-algebra
which is \'etale of rank~three. Further,~let
$l := c_1 T_1 + \ldots + c_4 T_4$
be a linear form with coefficients
in~$A$
and~$u \in D$.

Then,~we say that the cubic form
$\NT_u(l) := \tr (u \!\cdot\! \N (l))$
is obtained
from~$l$
and~$u$
by the {\em norm-trace~construction.} Correspondingly~for the cubic
surface~$S_u(l)$
over~$\bbQ$
given
by~$\NT_u(l) = 0$.
\end{defio}

\begin{rems}
\begin{iii}
\item
Actually,~$D$
is either a quadratic number field or isomorphic
to~$\bbQ \oplus \bbQ$.
In~the latter case, we simply have two \'etale
$\bbQ$-algebras
$A_0$
and~$A_1$,
both of rank~three. We~start with two linear forms
$l_0$
and~$l_1$
with coefficients
in~$A_0$
and~$A_1$,~respectively.
The~norm-trace~construction then degenerates to
$$u_0 \N(l_0) + u_1 \N(l_1)$$
for~$u = (u_0, u_1)$.
\item
An~\'etale algebra of rank three over a
field~$K$
may be

$\bullet$
$K \times K \times K$,

$\bullet$
the direct product of
$K$
with a quadratic field~extension, or

$\bullet$
a cubic field extension.

Further,~the cubic field extension may be Galois or~not. In~other words, the corresponding Galois group may be
$A_3$
or~$S_3$.
In~this language, the degenerate cases correspond to the non-transitive subgroups
$\bbZ/2\bbZ$
and~$0$
of~$S_3$.
\item
$A$
is actually always a free
$D$-module.
Indeed,~$D$
is a semilocal~ring.
Actually,~$D$
is even~Artin. Hence,~every locally free module of finite rank is~free.
\item
As~an \'etale
$\bbQ$-algebra
of rank two,
$D$~allows
two algebra homomorphisms
$\iota_0, \iota_1 \colon D \to \overline\bbQ$.

On~the other hand, as~a
$\bbQ$-algebra,
$A$~is
\'etale,~too. This~means,
$A$~is
a commutative semisimple~algebra of rank~six. There~are exactly six algebra homomorphisms
from~$A$
to~$\overline\bbQ$.
Three~of them are extensions
of~$\iota_0$,
the others
of~$\iota_1$.
We~denote them by
$\tau_0$,~$\tau_1$,~$\tau_2$
and
$\tau_3$,
$\tau_4$,
$\tau_5$,
respectively.

In~these terms, the norm-trace construction, applied
to~$l = c_1 T_1 + \ldots + c_4 T_4$
and~$u$,
yields the cubic~form
$\NT_u(l) = \iota_0(u) l^{\tau_0} l^{\tau_1} l^{\tau_2} + \iota_1(u) l^{\tau_3} l^{\tau_4} l^{\tau_5}$.
More~explicitly, this~is
\begin{eqnarray*}
&&
\iota_0(u) [\tau_0(c_1) T_1 + \!\ldots\! + \tau_0(c_4) T_4][\tau_1(c_1) T_1 + \!\ldots\! + \tau_1(c_4) T_4][\tau_2(c_1) T_1 + \!\ldots\! + \tau_2(c_4) T_4] \\
&&\hspace{-5mm} {} +
\iota_1(u) [\tau_3(c_1) T_1 + \!\ldots\! + \tau_3(c_4) T_4][\tau_4(c_1) T_1 + \!\ldots\! + \tau_4(c_4) T_4][\tau_5(c_1) T_1 + \!\ldots\! + \tau_5(c_4) T_4] \, .
\end{eqnarray*}
\end{iii}
\end{rems}

\begin{prop}
\label{trace}
Suppose,~we are given a commutative semi\-simple\/
\mbox{$\bbQ$-}al\-ge\-bra\/~$D$
of dimension~two and a commutative\/
$D$-algebra\/~$A$
which is \'etale of rank~three.
Further,~let\/
$u \in D$
and\/
$l$~be
a linear form in four
variables\/~$T_1, \ldots, T_4$
with coefficients
in\/~$A$.\smallskip

\noindent
Denote~by\/~$d$
the dimension of the
$\bbQ$-vector~space\/~$\langle l^{\tau_0}, \ldots, l^{\tau_5} \rangle \subseteq \Gamma(\bP^3, \calO(1))$.
Fix,~finally, a
field\/~$K$
such
that\/~$K \supseteq \im \tau_0, \ldots, \im \tau_5$.\smallskip\pagebreak[3]

\noindent
Then,

\begin{iii}
\item
$l^{\tau_0}, \ldots , l^{\tau_5}$~define
a rational~map
\begin{center}\bigskip
$\underline\iota \colon S_u(l) \times_{\Spec \bbQ} \Spec K \;\,\ratarrow\; \bP^5_{\!K} \, .$
\end{center}\bigskip
The~image
of\/~$\underline\iota$
is contained in a linear subspace of
dimension\/~$d-1$.
\item
If\/~$d = 4$
then\/
$S_u(l)$~is
a cubic surface
over\/~$\bbQ$
such that\/
$S_u(l) \times_{\Spec \bbQ} \Spec K$
is in generalized Cayley-Salmon~form.

More~precisely, if\/
$a_0 l^{\tau_0} + \ldots + a_5 l^{\tau_5} = 0$
and\/~$b_0 l^{\tau_0} + \ldots + b_5 l^{\tau_5} = 0$
are linearly independent relations then\/
$\underline\iota$
induces an~isomorphism
$$\iota \colon S_u(l) \times_{\Spec \bbQ} \Spec K \stackrel{\cong}{\longrightarrow} S_{\iota_0(u), \iota_1(u)}^{(a_0, \ldots, a_5, b_0, \ldots, b_5)} \, .$$
\item
If\/~$d \leq 3$
then\/
$S_u(l)$~is
the cone over a, possibly degenerate, cubic~curve.
\end{iii}\smallskip

\noindent
{\bf Proof.}
{\em
i)~is~standard.\smallskip

\noindent
ii)
In~this case, the forms
$l^{\tau_0}, \ldots, l^{\tau_5}$
generate the
\mbox{$K$-vector}
space
$\Gamma (\bP^3_{\!K}, \calO(1))$
of all linear~forms. Therefore,~they define a closed
immersion
of~$\bP^3_{\!K}$
into~$\bP^5_{\!K}$.
In~particular,
$\iota$~is
a closed~immersion.

We~have the relations
$a_0 l^{\tau_0} + \ldots + a_5 l^{\tau_5} = 0$
and~$b_0 l^{\tau_0} + \ldots + b_5 l^{\tau_5} = 0$.
The~cubic
surface~$S_u(l) \times_{\Spec \bbQ} \Spec K \subset \bP^3_{\!K}$
is given by
$\iota_0(u) l^{\tau_0} l^{\tau_1} l^{\tau_2} + \iota_1(u) l^{\tau_3} l^{\tau_4} l^{\tau_5} = 0$.
Consequently,~$\underline\iota$
maps
$S_u(l) \times_{\Spec \bbQ} \Spec K$~to
the cubic surface in generalized Cayley-Salmon~form
$\smash{S_{\iota_0(u), \iota_1(u)}^{(a_0, \ldots, a_5, b_0, \ldots, b_5)} \subset \bP^5_{\!K}}$.\smallskip

\noindent
iii) is~clear.
}
\eop
\end{prop}

\begin{ttt}
As~an application of the norm-trace construction, we have an explicit version of Galois~descent. For~this, some notation has to be~fixed.
\end{ttt}

\begin{nota}
For~$\sigma \in \Gal(\overline\bbQ/\bbQ)$,
denote 
by~$t_\sigma \colon \Spec \overline\bbQ \to \Spec \overline\bbQ$
the morphism of schemes induced
by~$\sigma^{-1} \colon \overline\bbQ \leftarrow \overline\bbQ$.
This~yields a morphism
$$t_\sigma^{\bP^5} \colon \bP^5_{\overline\bbQ} \longrightarrow \bP^5_{\overline\bbQ}$$
of
\mbox{$\overline\bbQ$-schemes}
which is
{\em twisted
by\/}~$\sigma$.
I.e.,~compatible
with~$\smash{t_\sigma \colon \Spec \overline\bbQ \to \Spec \overline\bbQ}$.
Observe~that,
on~$\overline\bbQ$-rational~points,
$$t_\sigma^{\bP^5} \colon (x_0 : \ldots : x_5) \mapsto \big( \sigma(x_0) : \ldots : \sigma(x_5) \big) \, .$$
We~will usually
write~$t_\sigma$
instead
of~$t_\sigma^{\bP^5}$.
The~morphism~$t_\sigma$
maps the cubic
surface~$S_{u_0, u_1}^{(a_0, \ldots, a_5, b_0, \ldots, b_5)}$
to~$S_{\sigma(u_0), \sigma(u_1)}^{(\sigma(a_0), \ldots, \sigma(a_5), \sigma(b_0), \ldots, \sigma(b_5))}$.\smallskip

Suppose,~for an \'etale
algebra~$A$
of rank~three over a commutative semi\-simple\/
\mbox{$\bbQ$-}al\-ge\-bra\/~$D$
of dimension~two and elements
$u \in D$
and~$a, b \in A$,
we have
$u_0 = \iota_0(u)$,
$u_1 = \iota_1(u)$,
$$a_i = \tau_i (a), \quad {\rm and} \quad b_i = \tau_i (b)$$
for~$i = 0, \ldots, 5$.
Assume~that
$a_0, \ldots, a_5$
are pairwise different from each~other.
Then,~every
$\sigma \in \Gal(\overline\bbQ/\bbQ)$
uniquely determines a permutation
$\pi_\sigma \in S_6$
such~that
$$\sigma(a_i) = a_{\pi_\sigma (i)} \quad {\rm and} \quad \sigma(b_i) = b_{\pi_\sigma (i)} \, .$$
This~yields a group
homomorphism~$\Pi \colon \Gal(\overline\bbQ/\bbQ) \to S_6$.
We~will denote the automorphism
of~$\bP^5$,
given by the
permutation~$\pi$
on coordinates,
by~$\pi$,~too.

Observe~that the permutations
$\pi_\sigma \in S_6$
preserve the block
structure~$[123][456]$.
Indeed,~$\sigma$
may eigher interchange the two algebra homomorphisms 
$\iota_0, \iota_1: D \to \bbC$
or~not. As~a consequence of this, we see~that
$\pi_\sigma$
defines a
morphism~$S_{\sigma(u_0), \sigma(u_1)}^{(\sigma(a_0), \ldots, \sigma(a_5), \sigma(b_0), \ldots, \sigma(b_5))} \to S_{u_0, u_1}^{(a_0, \ldots, a_5, b_0, \ldots, b_5)}$.

Putting~everything together,
$$\pi_\sigma \!\circ\! t_\sigma \colon S_{u_0, u_1}^{(a_0, \ldots, a_5, b_0, \ldots, b_5)} \longrightarrow S_{u_0, u_1}^{(a_0, \ldots, a_5, b_0, \ldots, b_5)}$$
is an automorphism
twisted
by~$\sigma$.
These~automorphisms form an operation
of~$\Gal(\overline\bbQ/\bbQ)$
on~$S_{u_0, u_1}^{(a_0, \ldots, a_5, b_0, \ldots, b_5)}$
from the~left.
\end{nota}

\begin{theoo}[{\rm Explicit Galois descent}{}]
\label{descent}
~\\[0pt]
Let\/~$A$
be an \'etale
algebra of rank~three over a commutative semi\-simple\/
\mbox{$\bbQ$-}al\-ge\-bra\/~$D$
of dimension~two. For~elements\/
$u \in D$
and\/
$a, b \in A$,
put\/
$u_0 := \iota_0(u)$,
$u_1 := \iota_1(u)$,
and,
for\/~$i = 0, \ldots, 5$,
$$a_i = \tau_i (a) \quad {\rm and} \quad b_i = \tau_i (b) \, .$$
Suppose~that\/
$a$
and\/
$b$
are\/
$\bbQ$-linearly~independent.
Assume~further that\/
$a_0, \ldots, a_5$
are pairwise different from each~other.

\begin{iii}
\item
Then,~there exist a cubic
surface\/~$S = S^{u_0, u_1}_{(a_0, \ldots, a_5, b_0, \ldots, b_5)}$
over\/~$\bbQ$
and an~isomorphism
$$\iota \colon S \times_{\Spec \bbQ} \Spec \overline\bbQ \stackrel{\cong}{\longrightarrow} S_{u_0, u_1}^{(a_0, \ldots, a_5, b_0, \ldots, b_5)}$$
such that, for every\/
$\sigma \in \Gal(\overline\bbQ/\bbQ)$,
the diagram
$$
\diagram
S \times_{\Spec \bbQ} \Spec \overline\bbQ \rrto^{\;\;\;\;\;\;\;\;\;\;\iota} & & S_{u_0, u_1}^{(a_0, \ldots, a_5, b_0, \ldots, b_5)} \\
S \times_{\Spec \bbQ} \Spec \overline\bbQ \rrto^{\;\;\;\;\;\;\;\;\;\;\iota} \uto^{\id \times t_\sigma} & & S_{u_0, u_1}^{(a_0, \ldots, a_5, b_0, \ldots, b_5)} \uto_{\pi_\sigma \circ t_\sigma}
\enddiagram
$$
commutes.
\item
The~properties given
determine\/~$S$
up to a unique isomorphism
of\/~\mbox{$\bbQ$-schemes}.
\item
Explicitly,~the\/
\mbox{$\bbQ$-scheme\/}~$S$
may be obtained by the norm-trace~construction as~follows.
$$S := S_u(l)$$
for\/~$l = c_1 T_1 + \ldots + c_4 T_4$
any linear form such that\/
$\tr (al) = 0$,
$\tr (bl) = 0$,
and\/
$c_0, \ldots, c_3 \in A$
are linearly independent
over\/~$\bbQ$.
\end{iii}
\smallskip

\noindent
{\bf Proof.}
{\em
i) and~ii)
These~assertions are particular cases of standard results from the theory of Galois~descent~\cite[Chapitre~V, \S4, n$^\circ$\,20, or~J,~Proposition~2.5]{Se}. In~fact, the scheme
$\smash{S_{u_0, u_1}^{(a_0, \ldots, a_5, b_0, \ldots, b_5)}}$
is (quasi-)projective
over~$\overline\bbQ$.
Thus,~everything which is needed are ``descent~data'', an operation
of~$\Gal(\overline\bbQ/\bbQ)$
on~$\smash{S_{u_0, u_1}^{(a_0, \ldots, a_5, b_0, \ldots, b_5)}}$
such that
$\sigma \in \Gal(\overline\bbQ/\bbQ)$
acts by a morphism
of~$\overline\bbQ$-schemes
which is twisted
by~$\sigma$.\smallskip

\noindent
iii)
The~$\bbQ$-linear
system of equations
\begin{eqnarray*}
\tr(ac) & = & 0 \, , \\
\tr(bc) & = & 0
\end{eqnarray*}
has a four-dimensional
space~$\bbL$
of~solutions. Indeed,~the bilinear form
$(x, y) \mapsto \tr(xy)$
is non-degenerate~\cite[\S8, Proposition~1]{Bou}. Hence,~the first two conditions
on~$l$
express
that~$c_0, \ldots , c_3 \in \bbL$
while the last one is equivalent to saying that
$\langle c_0, \ldots , c_3 \rangle$~is
a basis of that~space.

To~exclude the possibility that
$S$~degenerates
to a cone and to obtain the
isomorphism~$\iota$,
we intend to use Proposition~\ref{trace}.ii). This~requires to show that the linear forms
$l^{\tau_i} = c_1^{\tau_i} T_1 + \ldots + c_4^{\tau_i} T_4$
for~$0 \leq i \leq 5$
form a generating system of the vector space of all linear~forms. Equivalently,~we claim that the
\mbox{$6 \!\times\! 4$-matrix}
$$(c_j^{\tau_i})_{0 \leq i \leq 5, 1 \leq j \leq 4}$$
is of
rank~$4$.

To~prove this, we extend
$\{ c_1, \ldots, c_4 \}$
to a \mbox{$\bbQ$-basis}
$\{ c_1, \ldots, c_6 \}$
of~$A$.
It~is enough to verify that the
\mbox{$6 \!\times\! 6$-matrix}
$\smash{(c_j^{\tau_i})_{0 \leq i \leq 5, 1\leq j \leq 6}}$
is of full~rank. This~assertion is actually independent of the particular choice of a~basis. We~may do the calculations as well
with~$\smash{\{ 1, a, \ldots, a^5 \}}$.
This~yields the Vandermonde~matrix
$$
\left(
\begin{array}{cccc}
1 & a^{\tau_0} & \cdots & (a^{\tau_0})^5 \\[-2pt]
\vdots &\vdots   & \ddots & \vdots           \\[2pt]
1 & a^{\tau_5} & \cdots & (a^{\tau_1})^5 
\end{array}
\right) .
$$
of determinant equal~to
$$\smash{\prod_{i < j}} (a^{\tau_i} - a^{\tau_j}) = \smash{\prod_{i < j}} (a_i - a_j) \neq 0 \, .$$
Observe~that the six algebra homomorphisms
$\tau_i \colon A \to \overline\bbQ$
are given
by~$a \mapsto a_i$.

Consequently,~the linear
forms~$l^{\tau_i}$
yield the desired isomorphism
$$\iota \colon S \times_{\Spec \bbQ} \Spec K \stackrel{\cong}{\longrightarrow} S_{u_0, u_1}^{(a_0, \ldots, a_5, b_0, \ldots, b_5)} \, .$$
Indeed,~we have the
equations~$\tr (a c_i) = 0$.
Explicitly,~they express that, for each
$i \in \{1, \ldots, 4\}$,
$$0 = (a c_i)^{\tau_0} + \ldots + (a c_i)^{\tau_5} = a_0 c_i^{\tau_0} + \ldots + a_5 c_i^{\tau_5} \, .$$
This~means~$a_0 l^{\tau_0} + \ldots + a_5 l^{\tau_5} = 0$.
Analogously,~we
have~$b_0 l^{\tau_0} + \ldots + b_5 l^{\tau_5} = 0$.\smallskip

It~remains to verify the commutativity of the~diagram. For~this, we cover
$S_{u_0, u_1}^{(a_0, \ldots, a_5, b_0, \ldots, b_5)}$
by the affine open subsets given
by~$X_j \neq 0$
for~$j = 0, \ldots, 5$.
Observe~that the morphisms to be compared are both morphisms
of~\mbox{$\overline\bbQ$-schemes}
twisted
by~$\sigma$.
Hence,~we may compare the pull-back maps between the algebras of regular functions by testing their~generators.

For~arbitrary~$i \neq j$,
consider the rational
function~$X_i/X_j$.
Its~pull-back
under~$\iota$
is~$l^{\tau_i} / l^{\tau_j}$.
Therefore,~the pull-back
of~$X_i/X_j$
along the upper left corner~is
$$l^{\sigma^{-1} \circ \tau_i} / l^{\sigma^{-1} \circ \tau_j} = (c_1^{\sigma^{-1} \circ \tau_i} t_1 + \ldots + c_4^{\sigma^{-1} \circ \tau_i} t_4) / (c_1^{\sigma^{-1} \circ \tau_j} t_1 + \ldots + c_4^{\sigma^{-1} \circ \tau_j} t_4) \, .$$
On~the other hand, the pull-back
of~$X_i/X_j$
under~$\pi_\sigma \!\circ\! t_\sigma$
is~$X_{\pi_{\sigma^{-1}} (i)} / X_{\pi_{\sigma^{-1}} (j)}$.
Consequently,~for the pull-back along the lower right corner, we~find
\begin{eqnarray*}
l^{\tau_{\pi_{\sigma^{-1}} (i)}} / l^{\tau_{\pi_{\sigma^{-1}} (j)}} & = & (c_1^{\tau_{\pi_{\sigma^{-1}} (i)}} t_1 + \ldots + c_4^{\tau_{\pi_{\sigma^{-1}} (i)}} t_4) / (c_1^{\tau_{\pi_{\sigma^{-1}} (j)}} t_1 + \ldots + c_4^{\tau_{\pi_{\sigma^{-1}} (j)}} t_4) \, , \\
 & = & (c_1^{\sigma^{-1} \circ \tau_i} t_1 + \ldots + c_4^{\sigma^{-1} \circ \tau_i} t_4) / (c_1^{\sigma^{-1} \circ \tau_j} t_1 + \ldots + c_4^{\sigma^{-1} \circ \tau_j} t_4) \, .
\end{eqnarray*}
Indeed,~the embeddings
$\tau_{\pi_{\sigma^{-1}} (i)}, \sigma^{-1} \!\circ\! \tau_i \colon A \to \overline\bbQ$
are the same as one may check on the
generator~$T$,
$$\tau_{\pi_{\sigma^{-1}} (i)} (T) = a_{\pi_{\sigma^{-1}} (i)} = \sigma^{-1} (a_i) = \sigma^{-1} (\tau_i (T)) = (\sigma^{-1} \!\circ\! \tau_i) (T) \, .$$
This~completes the~proof.
}
\eop
\end{theoo}

\begin{algoo}[{\rm Computation of the Galois descent}{}]
\label{expl}
Let~two polynomials
$g \in \bbQ[U]$
and
$f \in D[V]$
of degrees two and three, respectively, be given which define \'etale algebras
$D := \bbQ[U]/(g)$
and~$A := D[V]/(f)$.
Further,~let
$u \in D$
and
$a,b \in A$
be as in~Theorem~\ref{descent}.\smallskip

\noindent
Then,~this algorithm computes the Galois descent of the cubic surface
$S_{u_0, u_1}^{(a_0, \ldots, a_5, b_0, \ldots, b_5)}$
for
$u_0 = \iota_0 (u)$,
$u_1 = \iota_1 (u)$,
$a_i = \tau_i (a)$
and~$b_i = \tau_i (b)$.

\begin{iii}
\item
Compute,~according to the definition, the traces
$a_{ij} := \tr (aU^iV^j)$
and
$b_{ij} := \tr (bU^iV^j)$
for
$i = 0,1$
and~$j = 0,1,2$.
\item
Determine~the kernel of of the
$2 \!\times\! 6$-matrix
$$\left(
\begin{array}{cccccc}
a_{00} & a_{01} & a_{02} & a_{10} & a_{11} & a_{12} \\
b_{00} & b_{01} & b_{02} & b_{10} & b_{11} & b_{12}
\end{array}
\right) .$$
Choose~linearly independent kernel vectors
$(c^{00}_k, c^{01}_k, c^{02}_k, c^{10}_k, c^{11}_k, c^{12}_k) \in \bbQ^6$
for
$k = 1, \ldots, 4$.
\item
Compute the norm
$$\N_{A[T_1,\ldots,T_4]/D[T_1,\ldots,T_4]} \bigg[ \sum_{k=1}^4 \Big( \!\!\!\sum_{i=0,1\atop j=0,1,2} \!\!\!c^{ij}_k U^i V^j \Big) T_k \bigg] \, .$$
\item
Multiply~that cubic form with coefficients
in~$D$
by~$u$.
\item
Finally,~apply the trace coefficient-wise and output the resulting cubic form
in~$T_1, \ldots, T_4$
with 20 rational~coefficients.
\end{iii}
\end{algoo}

\begin{rems}
\begin{abc}
\item
To~compute the norm of a
polynomial~$F$
with coefficients
in~$A$,
we treat
$A[T_1,\ldots,T_4]$
as a free
$D[T_1,\ldots,T_4]$-module
with basis
$1,V,V^2$.
We~establish the
$3 \times 3$-matrix
associated with the multiplication
by~$F$
map and compute its~determinant.

More~generally, observe that all the computations in steps i), iii), and~iv) are executed in the algebra
$A$
which is of dimension~six
over~$\bbQ$.
In~order to perform Algorithm~\ref{expl}, it is not necessary to realize the Galois~hull or any other large algebra on the~machine.
\item
The case that
$D = \bbQ \oplus \bbQ$
and~$A = \bbQ[V]/(f_1) \oplus \bbQ[V]/(f_2)$
is included by taking
$g(U) := U^2 - 1$~and
$$\textstyle f := \frac{1+U}2 f_1 + \frac{1-U}2 f_2 \, .$$
\end{abc}
\end{rems}

\section{The Galois operation on the descent variety}

\begin{lem}
\label{hilfsp}
Let\/~$A$
be an \'etale
algebra of rank~three over a commutative semi\-simple\/
\mbox{$\bbQ$-}al\-ge\-bra\/~$D$
of dimension~two and\/
$u \in D$
as well as\/
$a, b \in A$
as in Theorem~\ref{descent}.

\begin{abc}
\item
If\/~$D \cong \bbQ \oplus \bbQ$
then\/~$\Phi_{u_0, u_1}^{(a_0, \ldots, a_5, b_0, \ldots, b_5)} \in \bbQ[T]$.
\item
Otherwise,~if\/
$D \cong \bbQ(\sqrt{d})$
then\/~$\sqrt{d}\,\Phi_{u_0, u_1}^{(a_0, \ldots, a_5, b_0, \ldots, b_5)} \in \bbQ[T]$.
\end{abc}
\smallskip

\noindent
{\bf Proof.}
{\em
We~have
\begin{eqnarray*}
\Phi_{u_0, u_1}^{(a_0, \ldots, a_5, b_0, \ldots, b_5)} (T) & = & \frac1{\iota_0 (u)} (\tau_0 (a) + \tau_0 (b) T)(\tau_1 (a) + \tau_1 (b) T)(\tau_2 (a) + \tau_2 (b) T) \\
 & & \hspace{0.95cm} {} - \frac1{\iota_1 (u)} (\tau_3 (a) + \tau_3 (b) T)(\tau_4 (a) + \tau_4 (b) T)(\tau_5 (a) + \tau_5 (b) T) \, .
\end{eqnarray*}
The~homomorphisms~$\tau_i$
are acted upon
by~$\Gal(\overline\bbQ/\bbQ)$.
An~element~$\sigma \in \Gal(\overline\bbQ/\bbQ)$
may either interchange the homomorphisms
$\iota_0, \iota_1 \colon D \to \bbC$
or~not.
Correspondingly,~$\sigma$
interchanges the two blocks
${\tau_0, \tau_1, \tau_2}$
and~${\tau_3, \tau_4, \tau_5}$
or~not. Anyway,~it respects the block
structure~$[\tau_0 \tau_1 \tau_2][\tau_3 \tau_4 \tau_5]$.\smallskip

\noindent
a)
In~this case, no automorphism
of~$\overline\bbQ$
may interchange
$\iota_0$
and~$\iota_1$.
Hence,~$\Phi_{u_0, u_1}^{(a_0, \ldots, a_5, b_0, \ldots, b_5)}$
is
$\Gal(\overline\bbQ/\bbQ)$~invariant.\smallskip

\noindent
b)
Here,~the same argument shows that
$\Phi^{(a_0, \ldots, a_5, b_0, \ldots, b_5)}$
is
$\Gal(\overline\bbQ/\bbQ(\sqrt{d}))$~invariant.
Hence,~$\Phi_{u_0, u_1}^{(a_0, \ldots, a_5, b_0, \ldots, b_5)} \in \bbQ(\sqrt{d})[T]$.
As~the polynomial is anti-invariant under the conjugation
of~$\bbQ(\sqrt{d})$,
the assertion~follows.
}
\eop
\end{lem}

\begin{nota}
We~denote
by~$\Gamma$
the subgroup
of~$S_6$
of all elements preserving the block
structure~$[012][345]$.
$\Gamma$~is
of
order~$72$.
As~an abstract group,
$\Gamma \cong (S_3 \times S_3) \rtimes \bbZ/2\bbZ$.
\end{nota}

\begin{defi}
$\Gamma$~operates
on the set of all bijections
$\rho \colon \{0,1,2\} \to \{3,4,5\}$
as~follows.
For~$\sigma \in \Gamma$,
consider the conjugation
$\sigma \!\circ\! \rho \!\circ\! \sigma^{-1}$.
If~$\sigma$
interchanges the two blocks then invert this~bijection. Otherwise,~put
$i^\sigma := \sigma \!\circ\! \rho \!\circ\! \sigma^{-1}$.
\end{defi}

\begin{rems}
\begin{iii}
\item
$\rho \mapsto \rho^\sigma$~is
an operation
of~$\Gamma$
from the~left.
\item
If~the bijection
$\rho$
connects
$i$
with~$j$
then
$\rho^\sigma$
connects
$\sigma(i)$
with~$\sigma(j)$.
\end{iii}
\end{rems}

\begin{prop}
Let\/~$A$
be an \'etale
algebra of rank~three over a commutative semi\-simple\/
\mbox{$\bbQ$-}al\-ge\-bra\/~$D$
of dimension~two and\/
$u \in D$
as well as\/
$a, b \in A$
as in Theorem~\ref{descent}. Further,~suppose that\/
$S_{u_0, u_1}^{(a_0, \ldots, a_5, b_0, \ldots, b_5)}$
is non-singular.

\begin{abc}
\item
Then,~$\disc \Phi_{u_0, u_1}^{(a_0, \ldots, a_5, b_0, \ldots, b_5)} \in \bbQ^*$.
\item
Further,~on the descent
variety\/~$\smash{S^{u_0, u_1}_{(a_0, \ldots, a_5, b_0, \ldots, b_5)}}$
over\/~$\bbQ$,
there~are
\begin{iii}
\item
nine~obvious lines given~by
$$L_{\{i,j\}} \colon \iota^* X_i = \iota^* X_j = 0$$
for\/
$i = 0,1,2$
and\/~$j = 3,4,5$,
\item
18~non-obvious lines given~by
$$L^\lambda_\rho \colon \iota^* Z_0 + \iota^* Z_{\rho(0)} = \iota^* Z_1 + \iota^* Z_{\rho(1)} = 0$$
for\/
$\lambda$
a zero of the auxiliary polynomial\/
$\Phi_{u_0, u_1}^{(a_0, \ldots, a_5, b_0, \ldots, b_5)}$
and\/~$\rho \colon \{0,1,2\} \to \{3,4,5\}$
a~bijection. Here,~the
coordinates\/~$Z_i$
are given~by
\begin{eqnarray*}
Z_0 := -Y_0 + Y_1 + Y_2, \quad
Z_1 := Y_0 - Y_1 + Y_2, \quad
Z_2 := Y_0 + Y_1 - Y_2, \\
Z_3 := -Y_3 + Y_4 + Y_5, \quad
Z_4 := Y_3 - Y_4 + Y_5, \quad
Z_5 := Y_3 + Y_4 - Y_5 \hspace{1.5mm}
\end{eqnarray*}
for
$$Y_i := (a_i + b_i \lambda) X_i, \quad i = 0,\ldots,5\,.$$
\end{iii}
\item
An~element\/
$\sigma \in \Gal(\overline\bbQ/\bbQ)$
acts on the lines according to the~rules
$$\sigma(L_{\{i,j\}}) = L_{\{\pi_\sigma(i), \pi_\sigma(j)\}}, \qquad \sigma(L^\lambda_\rho) = L^{\lambda^\sigma}_{\rho^{\pi_\sigma}} \, .$$

\end{abc}\smallskip

\noindent
{\bf Proof.}
{\em
a)
By~Lemma~\ref{hilfsp}, the polynomial
$\Phi_{u_0, u_1}^{(a_0, \ldots, a_5, b_0, \ldots, b_5)}$
has rational coefficients after multiplication
with~$\sqrt{d}$
for a
suitable~$d \in \bbQ^*$.
This~implies
$\disc \Phi_{u_0, u_1}^{(a_0, \ldots, a_5, b_0, \ldots, b_5)} \in \bbQ$
as the cubic discriminant is homogeneous of degree four in the~coefficients. According~to Proposition~\ref{34}.ii), smoothness implies
$\smash{\disc \Phi_{u_0, u_1}^{(a_0, \ldots, a_5, b_0, \ldots, b_5)} \neq 0}$.\smallskip

\noindent
b)
The~isomorphism
$$\iota\colon S^{u_0, u_1}_{(a_0, \ldots, a_5, b_0, \ldots, b_5)} \times_{\Spec \bbQ} \Spec \overline\bbQ \longrightarrow S_{u_0, u_1}^{(a_0, \ldots, a_5, b_0, \ldots, b_5)}$$
is~provided by~Theorem~\ref{descent}. We~therefore obtain all the lines by pull-back
from~$S_{u_0, u_1}^{(a_0, \ldots, a_5, b_0, \ldots, b_5)}$.
The~formulas for them are given in~\ref{41} and Proposition~\ref{43}.b). Observe~Proposition~\ref{44} which ensures that the lines given are~distinct.\smallskip

\noindent
c)
From~the commutative diagram given in Theorem~\ref{descent}.i), we see that the operation
of~$\sigma$
on~$S^{u_0, u_1}_{(a_0, \ldots, a_5, b_0, \ldots, b_5)} \times_{\Spec \bbQ} \Spec \overline\bbQ$
goes over into the~automorphism
$$\pi_\sigma \!\circ\! t_\sigma \colon S_{u_0, u_1}^{(a_0, \ldots, a_5, b_0, \ldots, b_5)} \longrightarrow S_{u_0, u_1}^{(a_0, \ldots, a_5, b_0, \ldots, b_5)} \, .$$
$\pi_\sigma$~permutes
the coordinates while
$t_\sigma$~is
the operation
of~$\sigma$
on the~coefficients.
}
\eop
\end{prop}

\begin{remark}
The~nine obvious lines are defined by the pair of Steiner~trihedra
$$\big( \{ \iota^* E_0,  \iota^* E_1, \iota^* E_2 \},
        \{ \iota^* E_3,  \iota^* E_4, \iota^* E_5 \} \big)$$
which is clearly Galois~invariant.

If~$D \cong \bbQ \oplus \bbQ$
then the Steiner trihedra are Galois invariant~themselves. Otherwise,~they are permuted by the conjugation of that quadratic number~field.
\end{remark}

\begin{coro}
Let\/~$A$
be an \'etale
algebra of rank~three over a commutative semi\-simple\/
\mbox{$\bbQ$-}al\-ge\-bra\/~$D$
of dimension~two and\/
$u \in D$
as well as\/
$a, b \in A$
as in Theorem~\ref{descent}. Further,~suppose that\/
$S_{u_0, u_1}^{(a_0, \ldots, a_5, b_0, \ldots, b_5)}$
is non-singular.\smallskip

\noindent
If,~in this situation, the auxiliary polynomial has a rational zero or degenerates to a quadratic polynomial then\/
$S^{u_0, u_1}_{(a_0, \ldots, a_5, b_0, \ldots, b_5)}$
has a Galois invariant double-six.\medskip

\noindent
{\bf Proof.}
{\em
First~note that
$\smash{\Phi_{u_0, u_1}^{(a_0, \ldots, a_5, b_0, \ldots, b_5)}}$
can not degenerate to a polynomial of degree less than two as we have
$\disc (\Phi_{u_0, u_1}^{(a_0, \ldots, a_5, b_0, \ldots, b_5)}) \neq 0$
in the smooth~case. If~it degenerates to a quadratic polynomial then we have six non-obvious lines corresponding to the
zero~$\lambda = \infty$.

The~assumption therefore implies in any case that there is a Galois invariant set consisting of the nine obvious and six non-obvious~lines. The~complement of that is a double-six.
}
\eop
\end{coro}

\begin{coro}
\label{gerade}
Let\/~$A = D[V]/(f)$
be an \'etale
algebra of rank~three over a commutative semi\-simple\/
\mbox{$\bbQ$-}al\-ge\-bra\/~$D$
of dimension~two and\/
$u \in D$,
$a, b \in A$
as in Theorem~\ref{descent}. Further,~suppose that\/
$S_{u_0, u_1}^{(a_0, \ldots, a_5, b_0, \ldots, b_5)}$
is non-singular.
Let\/~$K \subset \overline\bbQ$
be any~subfield.

\begin{iii}
\item
Then,~$\Gal(\overline\bbQ/K)$
operates on the 45~tritangent~planes by even permutations if and only~if
$$\sqrt{\disc(\Phi_{u_0, u_1}^{(a_0, \ldots, a_5, b_0, \ldots, b_5)}) \,\disc(D)} \in K \, .$$
\item
$\Gal(\overline\bbQ/K)$
does not interchange the two pairs of Steiner trihedra complementary to the distinguished one if and only~if
$$\sqrt{\N_{D/\bbQ} (\disc f)} \in K \, .$$
\end{iii}

\noindent
{\bf Proof.}
{\em
i) We~will show this result in several~steps.\medskip

\noindent
{\em First step.}
Assume~that
$\sigma \in \Gal(\overline\bbQ/\bbQ)$
neither interchanges the zeroes of the auxiliary polynomial nor the embeddings
$\iota_0$
and~$\iota_1$.
Then,~the permutation is~even.\medskip

\noindent
It~will suffice to verify this for the case
$\pi_\sigma = (01)$.
Then,~there are only three invariant lines, namely
$L_{\{2, j\}}$
for~$j = 3,4,5$.
It~is not hard to check that, for each
line~$L$
different from these three,
the equations defining
$L$
and~$\sigma(L)$
together form a system of rank~three. Hence,~we have twelve orbits each consisting of two lines with a point in~common.\smallskip

\noindent
According~to Fact~\ref{Schlaf}, the permutation induced on the tritangent~planes
is a product of
$16$~two-cycles
and, therefore,~even.\medskip

\noindent
{\em Second step.}
Assume~that
$\sigma \in \Gal(\overline\bbQ/\bbQ)$
interchanges two zeroes of the auxiliary polynomial but does not interchange
$\iota_0$
and~$\iota_1$.
Then,~the permutation is~odd.\medskip

\noindent
In~view of the first step, it suffices to verify this for the
case~$\pi_\sigma = \id$.
Then,~there are
$15$
invariant lines, the nine obvious ones and the six non-obvious ones corresponding to the invariant root of the auxiliary polynomial.\smallskip

\noindent
According~to Fact~\ref{Schlaf}, the permutation induced on the tritangent~planes
is a product of
$15$~two-cycles.
Hence,~it is~odd.\medskip

\noindent
{\em Third step.}
Assume~that
$\sigma \in \Gal(\overline\bbQ/\bbQ)$
stabilizes the zeroes of the auxiliary polynomial but interchanges
$\iota_0$
and~$\iota_1$.
Then,~the permutation is~odd.\medskip

\noindent
Without~restriction, assume
that~$\pi_\sigma = (03)(14)(25)$.
Then,~there are
$15$
invariant lines. These~are the three obvious ones
$L_{\{1,4\}}$,
$L_{\{2,5\}}$,
and~$L_{\{3,6\}}$
and four for each value
of~$\lambda$.
The~latter ones correspond to the
bijections~$\rho\colon \{0,1,2\} \to \{3,4,5\}$
which fulfill at least one of the three conditions
$0 \mapsto 3$,
$1 \mapsto 4$,
and~$2 \mapsto 5$.\smallskip

\noindent
Again,~according to Fact~\ref{Schlaf}, the permutation induced on the tritangent~planes
is a product of
$15$~two-cycles
and, therefore,~odd.\medskip\pagebreak[3]

\noindent
{\em Fourth step.}
Conclusion.\medskip

\noindent
We~see that the operation
of~$\sigma$
on the 45 tritangent planes is even if and only if those
on~$\{\iota_0, \iota_1\}$
and the three zeroes
of~$\Phi_{u_0, u_1}^{(a_0, \ldots, a_5, b_0, \ldots, b_5)}$
are of the same~parity. This~means exactly~that
$$\sqrt{\disc(\Phi_{u_0, u_1}^{(a_0, \ldots, a_5, b_0, \ldots, b_5)}) \,\disc(D)}$$
is~$\sigma$-invariant.\medskip

\noindent
ii)
Let~$\sigma \in \Gal(\overline\bbQ/\bbQ)$.
According~to Proposition~\ref{44}, the two complementary pairs of Steiner trihedra are stable if and only if
$\pi_\sigma$
does not map the bijections
$\{0,1,2\} \to \{3,4,5\}$
from
$e$
to~$o$.
If~$\pi_\sigma$
preserves the two blocks
$\{0,1,2\}$
and~$\{3,4,5\}$
then this means that
$\pi_\sigma$~must
be~even. Otherwise,~it must be~odd.
As~$A \cong \bbQ[V]/(\N_{D/\bbQ}f)$,
we have exactly that
$\smash{\sqrt{\disc(\N_{D/\bbQ}f) \,\disc(D)}}$
is~\mbox{$\sigma$-invariant}.\smallskip

\noindent
We~claim
that~$\disc(\N_{D/\bbQ}f) \,\disc(D)$
equals~$\N_{D/\bbQ} (\disc f)$
up to square~factors. Indeed,~the definition of the discriminant together with the definition of the resultant~\cite[p.~79]{CLO}, implies, for any two~polynomials,
$$\disc(fg) =  \disc(f) \disc(g) \Res(f,g)^2.$$
Hence,~$\disc(\N_{D/\bbQ}f) = \disc(f) \disc(\overline{f}) \Res(f,\overline{f})^2$.
But~$\Res(f,\overline{f})$
is anti-invariant under the conjugation
of~$D$.
}
\eop
\end{coro}

\begin{remark}
To~operate on the 45~tritangent~planes by even permutations is a property characterizing the index two subgroup
of~$W(E_6)$.
This~is the simple group of
order~$25\,920$.
On~the other hand, the operation on the 27~lines is always~even.
\end{remark}

\begin{fac}
\label{Schlaf}
In\/~$W(E_6)$,
there are exactly four conjugacy classes of elements of order~two. The~corresponding operations on the lines and tritangent planes of a non-singular cubic surface are as~follows.

\begin{iii}
\item
There~are\/
$15$
invariant~lines. Then,~there are\/
$15$
invariant tritangent~planes. There~are\/
$15$
further orbits which are~pairs.
\item
There~are seven invariant~lines. Then,~there are five invariant tritangent planes and\/
$20$
orbits which are formed by~pairs.
\item
There~are three invariant~lines. The~others form six orbits of two lines which are skew and six orbits of two lines with a point in~common.

In~this case, there are seven invariant tritangent planes and\/
$19$~pairs.
\item
There~are three invariant~lines. The~others form twelve orbits each consisting of two lines with a point in~common.

Then,~there are\/
$13$
invariant tritangent~planes. There~are\/
$16$
further orbits which are~pairs.
\end{iii}\medskip

\noindent
{\bf Proof.}
{\em
This~result is essentially due to L.~Schl\"af\/li~\cite[pp.~114f.]{Sch} in the guise of his classification of cubic surfaces over the real~field. Today,~the four conjugacy classes may be seen directly using~{\tt GAP}. The~operations on the lines and tritangent planes are easily described in the blown-up~model.
}
\eop
\end{fac}

\section{Examples}

\begin{ttt}
Using~Algorithm~\ref{expl}, we generated a series of examples of smooth cubic surfaces
over~$\bbQ$.
Our~list of examples realizes each of the 63 conjugacy classes of subgroups
of~$W(E_6)$
which fix a pair of Steiner~trihedra but no~double-six. It~is available at the web~page {\tt http:/www.uni-math.gwdg.de/jahnel} of the second~author. In~this section, we present a few cubic surfaces from our list which are, as we think, of particular~interest.
\end{ttt}

\begin{ttt}
To~be more precise, our strategy to generate an example for a particular
subgroup~$G \subseteq [(S_3 \times S_3) \rtimes \bbZ/2\bbZ] \times S_3$
was as~follows.

\begin{iii}
\item
If~the group stabilizes the two Steiner trihedra then work
with~$g(U) := U^2 - 1$,
i.e.~$D = \bbQ \oplus \bbQ$.
Otherwise,~fix a polynomial defining a quadratic number~field.
\item
Choose~a starting
polynomial~$f \in D[V]$
such that
\begin{abc}
\item[ $\bullet$ ]
$A := D[V]/(f)$
is \'etale of rank three
over~$D$,
\item[ $\bullet$ ]
the Galois group of the degree-six
polynomial~$\N_{D[V]/Q[V]} (f) \in \bbQ[V]$
is equal
to~$\pr_1 (G) \subseteq (S_3 \times S_3) \rtimes \bbZ/2\bbZ$.
\end{abc}
\item
Always~put~$b := 1$.
Further,~
let~$(u,a) \in D \times A$
run through all pairs up to certain~height. For~each such pair, calculate the auxiliary
polynomial~$\Phi$.
If~that has the property desired in order to yield exactly the
group~$G$
then compute the descent variety in explicit form and terminate~immediately.
\end{iii}
\end{ttt}

\begin{rems}
\begin{iii}
\item
The~desired property of the auxiliary polynomial might simply be that it be generic with Galois group
$S_3$
or~$A_3$.
Several~of the groups require more restrictive conditions, for example that
$\bbQ(\sqrt{\disc \Phi})$
coincide
with~$D$.
\item
We~feel that the actual starting polynomial, we used in the experiment, is somehow~irrelevant. Of~importance are the
algebra~$A$
it defines and the
element~$a \in A$.
As~in all our examples
$\bbQ[a] = A$,
we will present the algebras in the form
$A := D[V]/(f)$
for~$a = \overline{V}$.
\item
For~42 of the 63~groups, we could start with the
algebra~$D = \bbQ \oplus \bbQ$.
In~this case, actually
$f = (f_0, f_1)$
for~$f_0, f_1 \in \bbQ[V]$.
The~other 21~groups required a quadratic number~field.
\end{iii}
\end{rems}

\begin{ex}
Start~with
$g(U) := U^2 - 1$,
i.e.~$D := \bbQ \oplus \bbQ$,
$f_0(V) := V^3 + \frac12 V + 1$,
$f_1(V) := V^3 - 2V^2  + 5$,
$u_0 = 1$,
and~$u_1 = 2$.
Both~polynomials have Galois
group~$S_3$.
The resulting non-singular cubic
surface~$S$
is given by the~equation
\begin{eqnarray*}
 & & 18T_1^3 - 40T_1^2T_2 + 37T_1^2T_3 - 30T_1^2T_4 + 68T_1T_2^2 + 4T_1T_2T_3 - 36T_1T_2T_4 -
64T_1T_3^2 \\
 & & {} - 14T_1T_3T_4 + 38T_1T_4^2 - 24T_2^3 - 6T_2^2T_3 - 12T_2^2T_4 - 72T_2T_3^2 + 
64T_2T_3T_4 + 16T_2T_4^2 \\
 & & {} \hspace{7.6cm} + 31T_3^3 - 12T_3^2T_4 + 27T_3T_4^2 - 5T_4^3 = 0 \, .
\end{eqnarray*}
Here,~the auxiliary polynomial is generic with Galois
group~$S_3$.
Hence,~the Galois group operating on the 27~lines is the maximal
$S_3 \times S_3 \times S_3$
stabilizing the two Steiner trihedra of a~pair. We~have orbit
structure~$[9,18]$.
\end{ex}

\begin{ex}
Start~with
$g(U) := U^2 - 7$,
i.e.~$D := \bbQ (\sqrt{7})$,
and
$f(V) := V^3 + (1 - \sqrt{7})V^2 + (-1 + \sqrt{7})V + (5 - \sqrt{7})$.
This~is a polynomial with Galois
group~$S_3$.
Its~discriminant
is~$(810 \sqrt{7} - 2376)$,
a number of
norm~$1026^2$.
\mbox{Finally},~put~$u := \sqrt{7}$.
The~resulting non-singular cubic
surface~$S$
is given by the~equation
\begin{eqnarray*}
 & & -5T_1^3 + 5T_1^2T_2 + 5T_1^2T_3 + 3T_1T_2^2 - 5T_1T_2T_3 + 5T_1T_2T_4 \\
 & & {} \hspace{2cm} + 4T_1T_3^2 - T_1T_3T_4 - 6T_1T_4^2 - 6T_2^3 - 3T_2^2T_3 - 6T_2^2T_4 + 2T_2T_3^2 \\
 & & {} \hspace{4cm} + 2T_2T_3T_4 - 4T_2T_4^2 - 5T_3^3 - 4T_3^2T_4 - 4T_3T_4^2 - 2T_4^3 = 0 \, .
\end{eqnarray*}
The~auxiliary polynomial has Galois
group~$S_3$.

As~$N_{D/\bbQ}(\disc(f))$
is a perfect square, the Galois group operating on the 27~lines stabilizes three pairs of Steiner trihedra which are complementary in the sense that together they contain all the 27~lines. Actually,~it is the maximal group with this~property. It~is of index two
in~$[(S_3 \times S_3) \rtimes \bbZ/2\bbZ] \times S_3$.
We~have orbit
structure~$[9,9,9]$.
\end{ex}

\begin{ex}
\label{77}
Start~with
$g(U) := U^2 - 1$,
i.e.~$D := \bbQ \oplus \bbQ$,
$f_0(V) := V^3 + 3V^2 - 9V - 19$,
$f_1(V) := V^3 - 15V^2 + \frac{261}4 V - 85$,
$u_0 = 4$,
and~$u_1 = 1$.
Then,~Algorithm~\ref{descent} yields the non-singular cubic
surface~$S$
given by the~equation
\begin{eqnarray*}
 & & 9T_1^3 + 4T_1^2T_2 + 6T_1^2T_3 - 3T_1T_2^2 - 2T_1T_2T_3 - 3T_1T_3^2 - T_2^3 - 3T_2^2T_3 \\
 & & {} \hspace{3.3cm} - 3T_2^2T_4 - 6T_2T_3^2 + 2T_2T_3T_4 + 11T_2T_4^2 + T_3^3 - 3T_3T_4^2 - T_4^3 = 0 \, .
\end{eqnarray*}
Here,~the auxiliary polynomial has Galois
group~$A_3$.
Even~more, its splitting field coincides with that
of~$f_1$.
Hence,~the Galois group operating on the 27~lines is of order~three. We~have orbit
structure~$[3,3,3,3,3,3,3,3,3]$.
\end{ex}

\begin{remark}
This~group of order three is of particular~interest in connection with the Brauer
group~$\Br(S) = H^1 (\Gal(\overline\bbQ/\bbQ), \Pic(S))$
of~$S$.
This~is an important invariant in the arithmetic of cubic~surfaces. For~more information on its applications, the reader is referred to Yu.~I.~Manin's book~\cite{Ma}.

It~turns out that
$\Br(S)$
is completely determined by the group operating on the 27~lines. It~may take only five~\cite{SD,Co} values,
$0$,
$\bbZ/2\bbZ$,
$\bbZ/2\bbZ \times \bbZ/2\bbZ$,
$\bbZ/3\bbZ$,
and~$\bbZ/3\bbZ \times \bbZ/3\bbZ$.
A~calculation in~{\tt GAP} shows that
$\Br(S) = \bbZ/3\bbZ \times \bbZ/3\bbZ$
for exactly one of the 350 conjugacy classes of subgroups
in~$W(E_6)$.

This~group was discussed in~\cite{Ma}, already, as the group operating on the 27~lines
of~``$aT_1^3 + T_2^3 + T_3^3 + T_4^3 = 0$''
for~$a \in \bbQ(\zeta_3)$
a non-cube. Example~\ref{77} shows that it appears for cubic surfaces
over~$\bbQ$,~too.
\end{remark}

\begin{ex}
Start~with
$g(U) := U^2 - 2$,
i.e.~$D := \bbQ (\sqrt{2})$,
and
$f(V) := V^3 - \frac32\sqrt{2} V + \sqrt{2}$.
We~have
$F(V) := \N_{D[V]/Q[V]} (f) = V^6 - \frac92 V^2 + 6V - 2$,
the Galois group of which is
$[(S_3 \times S_3) \rtimes \bbZ/2\bbZ] \cap A_6$
of
order~$36$.
In~particular, the Galois group
of~$f$
itself
is~$A_3$.

It~turns out that,
for~$u := (-\sqrt{2} + 5)$,
the auxiliary polynomial has discriminant, up to square factors, equal
to~$2$.
Hence,~the Galois operation on the 45 tritangent planes is~even. The~Galois
group~$G$
is of
order~$108$.
On~the 27~lines, the orbit structure
is~$[9,18]$.
\end{ex}

\begin{remark}
This~example is of interest from the technical point of~view. Observe~the following~particularities.

The~$2$-Sylow~subgroup
of~$G$
is cyclic of order~four.
Actually,~$G$
is the maximal even subgroup
of~$W(E_6)$
stabilizing a pair of Steiner trihedra with this~property.

The~Galois group
of~$F(V)$
of
order~$36$
has the same
$2$-Sylow~subgroup,~already.
It~is generated by a permutation of the
form~$(0314)(25)$.
Further,~its
$3$-Sylow~subgroup
is~normal. Hence,~we may obtain the cyclic group of order~four as a~quotient. Consequently,~there must be a quadratic extension
of~$D = \bbQ(\sqrt{d})$
which is Galois and even cyclic
over~$\bbQ$.
This~causes limitations
on~$D$
due to the following~fact.
\end{remark}

\begin{fac}
Let\/~$L/\bbQ$
be a Galois extension which is cyclic of degree~four. Then,~the quadratic intermediate
field\/~$\bbQ(\sqrt{d})$
is real and the norm of the fundamental unit
is~$(-1)$.
In~particular,~$d \equiv 1,2 \pmod 4$.\smallskip

\noindent
{\bf Proof.}
{\em
We~have to show that
$(-1)$
is a norm
from~$D$.
For~this, according to the Hasse norm theorem, it is sufficient to verify that $(-1)$~is
a norm
from~$D_\nu$
for every
prime~$\nu$.
If~$[L_w : \bbQ_p] = 1$
then there is nothing to~prove.
If~$[L_w : \bbQ_p] = 2$
then the decomposition group
is~$2\bbZ/4\bbZ \subset \bbZ/4\bbZ$
which means that
$D_\nu = \bbQ_p$.
Thus,~there is nothing to prove,~either.
Finally,~if~$[L_w : \bbQ_p] = 4$
then, according to local class field theory,
$N L_w^* \subset N D_\nu^* \subset \bbQ_p^*$,
each of index two, such
that~$\bbQ_p^* / N L_w^* \cong \bbZ/4\bbZ$.
The~fact
that~$(-1)^2 = 1$
implies~$(-1) \in N D_\nu^*$.
}
\eop
\end{fac}

\end{document}